\documentclass[12pt,reqno]{amsart}

\usepackage{amsmath, epsfig, cite}
\usepackage{amssymb}
\usepackage{amsfonts}
\usepackage{latexsym}
\usepackage{graphicx}
\usepackage{algorithmic,algorithm}
\usepackage{amsthm,array}

\newtheorem{thm}{Theorem}[section]

\newtheorem{cor}[thm]{Corollary}

\newtheorem{lem}[thm]{Lemma}

\def\pf{\noindent {\it Proof.} }

\numberwithin{equation}{section}

\makeatletter \@addtoreset{equation}{section} \makeatother

\begin{document}
\allowdisplaybreaks
\title
{Congruences for the $k$ dots bracelet partition functions}
\author[Suping Cui]{Suping Cui}
\address[S. P. Cui]{Center for Combinatorics, LPMC-TJKLC, Nankai
University, Tianjin 300071, P. R. China} \email{jiayoucui@163.com}

\author[N. S. S. Gu]{Nancy Shanshan Gu}
\address[N. S. S. Gu]{Center for Combinatorics, LPMC-TJKLC, Nankai
University, Tianjin 300071, P. R. China} \email{gu@nankai.edu.cn}

\date{\today}
\begin{abstract} By finding the congruent relations between the generating function of the $5$ dots bracelet partitions
and that of the $5$-regular partitions, we get some new congruences
modulo $2$ for the $5$ dots bracelet partition function. Moreover,
for a given prime $p$, we study the arithmetic properties modulo $p$
of the $k$ dots bracelet partitions.
\end{abstract}

\maketitle

\section{Introduction}

In \cite{Andrews-Paule-2007}, Andrews and Paule studied the broken $k$-diamond partitions by using MacMahon's partition analysis,
and gave the generating function of $\Delta_k(n)$ which denotes the number of the broken $k$-diamond partitions of $n$:
\begin{equation}\label{k-diamond-gf}
\sum_{n=0}^{\infty}\Delta_k(n)q^n=\frac{(-q;q)_{\infty}}{(q;q)^2_{\infty}(-q^{2k+1};q^{2k+1})_{\infty}}.
\end{equation}
In \cite{Andrews-Paule-2007}, They proved the following arithmetic theorem for $\Delta_1(n)$.
\begin{thm}\label{delta-AP}\cite[Theorem 5]{Andrews-Paule-2007} For $n \geq 0$,
\begin{equation*}
\Delta_1(2n+1)\equiv 0 \quad (\text{mod}\ 3).
\end{equation*}
\end{thm}
Meanwhile, they posed some conjectures related to $\Delta_2(n)$. For other study of the arithmetic of the broken $k$-diamond partitions, see \cite{Chan-2008,Hirschhorn-Sellers-2007,Jameson-2013,Mortenson-2008,Paule-Radu-2010,Radu-Sellers-2011,Radu-Sellers-1-2012,Xiong}.
In \cite{Fu-2011}, Fu found a combinatorial proof of Theorem \ref{delta-AP} for $\Delta_1(n)$, and introduced a generalization of the broken
$k$-diamond partitions which he called the $k$ dots bracelet partitions. The number of this kind of partitions of $n$ is denoted by
$\mathfrak{B}_k(n)$, and the generating function of  $\mathfrak{B}_k(n)$ is stated as follows.
\begin{align*}
\sum^{\infty}_{n=0}\mathfrak{B}_{k}(n)q^n=\frac{(-q;q)_{\infty}}{(q;q)_{\infty}^{k-1}(-q^k;q^k)_{\infty}},\ k \geq 3.
\end{align*}

In \cite{Fu-2011}, Fu proved the following congruences for $\mathfrak{B}_k(n)$.
\begin{thm}\cite[Theorem 3.3]{Fu-2011} For $n>0$, $k\geq 3$, if $k=p^r$ is a prime power, we have
\begin{equation*}
\mathfrak{B}_k(2n+1) \equiv 0 \quad (\text{mod}\ p).
\end{equation*}
\end{thm}

\begin{thm}\cite[Theorem 3.5]{Fu-2011}\label{Fu-thm3.5} For any $k \geq 3$, $s$ an integer between $1$ and $p-1$ such that $12s+1$ is a quadratic nonresidue
modulo $p$, and any $n\geq 0$, if $p\mid k$ for some  prime $p \geq
5$  say $k=pm$, then we have
\begin{align*}
\mathfrak{B}_{k}(pn+s)&\equiv 0  \quad (\text{mod}\  p).
\end{align*}
\end{thm}

\begin{thm}\cite[Theorem 3.6]{Fu-2011} For $n\geq 0$, $k\geq 3$ even, say $k=2^{m}l$, where $l$ is odd, we have
\begin{align*}
\mathfrak{B}_{k}(2n+1)&\equiv 0  \ (\text{mod}\ 2^{m}).
\end{align*}
\end{thm}

Later, in \cite{Radu-Sellers-2012}, Radu and Sellers found some new congruences for $\mathfrak{B}_k(n)$.
\begin{thm}\cite[Theorem 1.4]{Radu-Sellers-2012} For all $n \geq 0$,
\begin{align*}
\mathfrak{B}_{5}(10n+7)&\equiv 0  \ (\text{mod}\  5^{2}),\\
\mathfrak{B}_{7}(14n+11)&\equiv 0  \ (\text{mod}\  7^{2}),\ \text{and}\\
\mathfrak{B}_{11}(22n+21)&\equiv 0  \ (\text{mod}\  11^{2}).
\end{align*}
\end{thm}

In this paper, we continue to study the arithmetic of the $k$ dots
bracelet partitions. First, we recall two kinds of partitions which
are used in this paper.

A partition of a positive integer $n$ is a nonincreasing sequence of positive
integers whose sum is $n$. Let $p(n)$ denote the number of partitions
of $n$. We know that
\begin{equation}\label{g-p}
\sum_{n=0}^{\infty}p(n)q^n=\frac{1}{(q;q)_{\infty}}.
\end{equation}

If $\ell$ is a positive integer, then a partition is called
$\ell$-regular partition if there is no part divisible by $\ell$.
Let $b_\ell(n)$ denote the number of $\ell$-regular partitions of
$n$. The generating function of $b_\ell(n)$ is stated as
follows.
\begin{equation}\label{g-l-regular}
\sum_{n=0}^{\infty}b_\ell(n)q^n=\frac{(q^\ell;q^\ell)_{\infty}}{(q;q)_{\infty}}.
\end{equation}

In section $2$, based on an identity given by Ramanujan in
\cite{Ramanujan-2000} and a congruence for the generating function
of $b_5(2n)$ given by Hirschhorn and Sellers in
\cite{Hirschhorn-Sellers-2010}, we obtain two congruences modulo $2$
for $B_5(n)$. Meanwhile, by finding a congruent relation between the
generating function of $\mathfrak{B}_5(n)$ and that of $b_5(n)$, we
get many infinite family of congruences modulo $2$ for
$\mathfrak{B}_5(n)$. In section $3$, for a given prime $p$, by means
of a $p$-dissection identity of $f(-q)$ given by the authors in
\cite{Cui-Gu-2012} and three classical congruences for $p(n)$ given
by Ramanujan in \cite{Ramanujan-1919,Ramanujan-2000}, we get more
congruences modulo $p$ for $\mathfrak{B}_k(n)$.

In the following, we list some definitions and identities which are
frequently used in this paper.

As usual, we follow the notation and terminology in
\cite{Gasper-Rahman-2004}. For $|q|<1$, the $q$-shifted factorial is
defined by
\begin{equation*}
(a;q)_\infty= \prod_{k=0}^{\infty}(1-aq^k) \quad\text{and}\quad
(a;q)_n = \frac{(a;q)_\infty}{(aq^n;q)_\infty}, \text{ for } n\in
\mathbb{C}.
\end{equation*}

The Legendre symbol is a function of of $a$ and $p$ defined as follows:
$$\left(\frac{a}{p}\right)=\left\{\begin{array}{ll}1,&\text{if } a \text{ is a quadratic residue modulo } p \text{ and } a \not\equiv 0\ (\text{mod}\ p),\\
-1,&\text{if } a \text{ is a quadratic non-residue modulo } p,\\
0,&\text{if } a\equiv 0\ (\text{mod}\ p). \end{array}\right.$$

Jacobi's triple product identity \cite[Theorem 1.3.3]{Berndt-2004}:
for $z \neq 0$ and $|q|<1$,
\begin{equation}\label{Jacobi}
\sum_{n=-\infty}^{\infty}z^nq^{n^2}=(-zq,-q/z,q^2;q^2)_{\infty}.
\end{equation}

Ramanujan's general theta function $f(a,b)$ is defined by
\begin{equation*}
f(a,b):=\sum_{n=-\infty}^{\infty}a^{\frac{n(n+1)}{2}}b^{\frac{n(n-1)}{2}},\qquad
|ab|<1.
\end{equation*}
A special case of $f(a,b)$ is stated as follows.
\begin{equation*}
f(-q):= f(-q,-q^2)=\sum_{n=-\infty}^{\infty}(-1)^nq^{\frac{n(3n-1)}{2}}=(q;q)_{\infty}.
\end{equation*}

\section{Congruences modulo $2$ for $\mathfrak{B}_5(n)$}

First, we recall an identity given by Ramanujan in \cite[p. 212]{Ramanujan-2000}.
\begin{equation}\label{Ramanujan-qq}
(q;q)_{\infty}=\frac{(q^{10},q^{15},q^{25};q^{25})_{\infty}}{(q^5,q^{20};q^{25})_{\infty}}-q(q^{25};q^{25})_{\infty}
-q^2\frac{(q^5,q^{20},q^{25};q^{25})_{\infty}}{(q^{10},q^{15};q^{25})_{\infty}}.
\end{equation}
For convenience, we set
\begin{equation*}
a(q)=\frac{(q^{10},q^{15};q^{25})_{\infty}}{(q^5,q^{20};q^{25})_{\infty}}\qquad \text{and} \qquad
b(q)=\frac{(q^5,q^{20};q^{25})_{\infty}}{(q^{10},q^{15};q^{25})_{\infty}}=\frac{1}{a(q)}.
\end{equation*}
Then, we can rewrite \eqref{Ramanujan-qq} as
\begin{equation}\label{Hir-qq}
(q;q)_{\infty}=(q^{25};q^{25})_{\infty}(a(q)-q-q^2b(q)).
\end{equation}

In \cite{Hirschhorn-Sellers-2010}, Hirschhorn and Sellers obtained the following congruence for the generating function of $b_5(2n)$.
\begin{equation}\label{b5555}
\sum_{n=0}^{\infty} b_{5}(2n)q^{n}\equiv(q^2;q^2)_\infty\quad (\text{mod}\ 2).
\end{equation}

By means of \eqref{Hir-qq} and \eqref{b5555}, we have the following results for $\mathfrak{B}_5(n)$.

\begin{thm} For $n\geq 0$, we have
\begin{align*}
\mathfrak{B}_5(10n+6)&\equiv 0 \quad (\text{mod}\ 2),\\
\mathfrak{B}_5(10n+8)&\equiv 0 \quad (\text{mod}\ 2).
\end{align*}
\end{thm}
\pf First, we have
\begin{align*}
\sum_{n=0}^\infty\mathfrak{B}_5(n)q^n&=\frac{(-q;q)_{\infty}}{(q;q)_{\infty}^4(-q^5;q^5)_{\infty}}\\
&=\frac{(q^2;q^2)_\infty(q^5;q^5)_\infty}{(q;q)^5_\infty(q^{10};q^{10})_\infty}\\
&\equiv \frac{(q^2;q^2)_\infty(q^5;q^5)_\infty}{(q^4;q^4)_\infty(q^{10};q^{10})_\infty(q;q)_\infty}\quad (\text{mod}\ 2)\\
&\equiv \frac{1}{(q^2;q^2)_\infty(q^{10};q^{10})_\infty}\frac{(q^5;q^5)_\infty}{(q;q)_\infty}\quad (\text{mod}\ 2)\\
&= \frac{1}{(q^2;q^2)_\infty(q^{10};q^{10})_\infty}\cdot \sum_{n=0}^\infty b_5(n)q^n.
\end{align*}
Then,
\begin{align*}
\sum_{n=0}^\infty\mathfrak{B}_5(2n)q^n&\equiv \frac{1}{(q;q)_\infty(q^{5};q^{5})_\infty}\cdot \sum_{n=0}^\infty b_5(2n)q^n\quad (\text{mod}\ 2)\\
&\equiv \frac{(q^2;q^2)_\infty}{(q;q)_\infty(q^{5};q^{5})_\infty}\quad (\text{mod}\ 2)\qquad \text{by } \eqref{b5555}\\
&\equiv \frac{(q;q)_\infty}{(q^{5};q^{5})_\infty}\quad (\text{mod}\ 2).
\end{align*}
According to \eqref{Hir-qq}, we have
\begin{equation}\label{B-Rama-5}
\sum_{n=0}^\infty\mathfrak{B}_5(2n)q^n \equiv \frac{(q^{25};q^{25})_{\infty}}{(q^{5};q^{5})_\infty}(a(q)-q-q^2b(q)) \quad (\text{mod}\ 2).
\end{equation}
Therefore, we get
\begin{align*}
\mathfrak{B}_5(2(5n+3))&= \mathfrak{B}_5(10n+6)\equiv0\quad (\text{mod}\ 2),\\
\mathfrak{B}_5(2(5n+4))&= \mathfrak{B}_5(10n+8)\equiv0\quad (\text{mod}\ 2).
\end{align*}\qed

\begin{lem}\label{B-b-2} We have
\begin{equation*}
\sum_{n=0}^\infty\mathfrak{B}_5(10n+2)q^n \equiv \sum_{n=0}^\infty b_5(n)q^n \quad (\text{mod}\ 2).
\end{equation*}
\end{lem}
\pf Due to \eqref{B-Rama-5}, we get
\begin{equation*}
\sum_{n=0}^\infty\mathfrak{B}_5(2(5n+1))q^n= \sum_{n=0}^\infty\mathfrak{B}_5(10n+2)q^n
\equiv\frac{(q^5;q^5)_\infty}{(q;q)_\infty}=\sum_{n=0}^\infty b_5(n)q^n\quad (\text{mod}\ 2).
\end{equation*}\qed

In \cite{Cui-Gu-2012}, the authors found many infinite family of congruences modulo $2$ for $b_5(n)$.
\begin{thm}\cite[Theorem 3.17]{Cui-Gu-2012}\label{cg-1} For any prime $p \geq 5$, $\left(\frac{-10}{p}\right)=-1$, $\alpha \geq 1$, and $n \geq 0$, we have
\begin{equation*}
b_5(4\cdot p^{2\alpha}n+\frac{(24i+7p)p^{2\alpha-1}-1}{6})\equiv 0 \quad (\text{mod}\ 2),
\end{equation*}
where $i=1,2,\ldots,p-1$.
\end{thm}
\begin{thm}\cite[Theorem 3.20]{Cui-Gu-2012}\label{cg-2} For $\alpha \geq 0$ and $n\geq 0$, we have
\begin{align*}
b_{5}(4\cdot5^{2\alpha+1}n+\frac{31\cdot 5^{2\alpha}-1}{6})&\equiv 0\quad (\text{mod}\ 2),\\
b_{5}(4\cdot5^{2\alpha+1}n+\frac{79\cdot 5^{2\alpha}-1}{6})&\equiv 0\quad (\text{mod}\ 2),\\
b_{5}(4\cdot5^{2\alpha+2}n+\frac{83\cdot 5^{2\alpha+1}-1}{6})&\equiv 0\quad (\text{mod}\ 2),\\
b_{5}(4\cdot5^{2\alpha+2}n+\frac{107\cdot 5^{2\alpha+1}-1}{6})&\equiv 0\quad (\text{mod}\ 2).
\end{align*}
\end{thm}
Therefore, combining Lemma \ref{B-b-2} with Theorem \ref{cg-1} and
Theorem \ref{cg-2}, we obtain some more congruences for
$\mathfrak{B}_5(n)$.
\begin{thm}\label{B5-p} For any prime $p \geq 5$, $\left(\frac{-10}{p}\right)=-1$, $\alpha \geq 1$, and $n \geq 0$, we have
\begin{equation*}
\mathfrak{B}_5(40\cdot p^{2\alpha}n+\frac{5\cdot (24i+7p)p^{2\alpha-1}+1}{3})\equiv 0 \quad (\text{mod}\ 2),
\end{equation*}
where $i=1,2,\ldots,p-1$.
\end{thm}
For example, by setting $p=17$, $i=6$, and $\alpha=1$ in Theorem \ref{B5-p}, we have the following congruence.
\begin{equation*}
\mathfrak{B}_5(11560n+7452)\equiv 0 \quad (\text{mod}\ 2).
\end{equation*}

\begin{thm}
For $\alpha \geq 1$ and $n\geq 0$, we have
\begin{align*}
\mathfrak{B}_5(8\cdot5^{2\alpha}n+\frac{31\cdot 5^{2\alpha-1}+1}{3})&\equiv 0\quad (\text{mod}\ 2),\\
\mathfrak{B}_5(8\cdot5^{2\alpha}n+\frac{79\cdot 5^{2\alpha-1}+1}{3})&\equiv 0\quad (\text{mod}\ 2),\\
\mathfrak{B}_5(8\cdot5^{2\alpha+1}n+\frac{83\cdot 5^{2\alpha}+1}{3})&\equiv 0\quad (\text{mod}\ 2),\\
\mathfrak{B}_5(8\cdot5^{2\alpha+1}n+\frac{107\cdot 5^{2\alpha}+1}{3})&\equiv 0\quad (\text{mod}\ 2).
\end{align*}
\end{thm}

\section{Congruences modulo $p$ for $\mathfrak{B}_k(n)$ }

In \cite{Cui-Gu-2012}, the authors studied a $p$-dissection identity of $f(-q)$ for a given prime $p\geq 5$.
\begin{thm}\cite[Theorem 2.2]{Cui-Gu-2012}\label{f-p} For any prime $p \geq 5$, we have
\begin{equation*}
f(-q)=\sum_{\tiny \begin{array}{l}k=-\frac{p-1}{2}\\k\neq \frac{\pm p-1}{6}\end{array}}^{\frac{p-1}{2}} (-1)^kq^{\frac{3k^{2}+k}{2}}f(-q^{\frac{3p^2+(6k+1)p}{2}},-q^{\frac{3p^2-(6k+1)p}{2}})+(-1)^{\frac{\pm p-1}{6}}q^{\frac{p^2-1}{24}}f(-q^{p^{2}}),
\end{equation*}
where $\pm$ depends on the condition that $(\pm p-1)/6$ should
be an integer. Meanwhile, we claim that $(3k^2+k)/2$ and
$(p^2-1)/24$ are not in the same residue class modulo $p$ for $-(p-1)/2
\leq k \leq (p-1)/2$ and $k \neq (\pm p-1)/6$.
\end{thm}
According to the above theorem, we have the following result.
\begin{lem}\label{lemma1} For any prime $p \geq 5,n\geq 0$, and $r\geq 1$, if $k=p^r$ is a prime power,
then for $1\leq \alpha \leq (r+1)/2$, we have
\begin{equation*}
\sum_{n=0}^{\infty}\mathfrak{B}_k(p^{2\alpha-1}n+\frac{p^{2\alpha}-1}{12})q^n \equiv \left((-1)^{\frac{\pm p-1}{6}}\right)^{\alpha}\frac{(q^{2p};q^{2p})_\infty}{(q^{2p^{r-(2\alpha-1)}};q^{2p^{r-(2\alpha-1)}})_\infty} \quad (\text{mod}\ p),
\end{equation*}
where $\pm$ depends on the condition that $(\pm p-1)/6$ should be an integer.
\end{lem}
\pf We prove the lemma by induction on $\alpha$.
For $k=p^r$, in \cite{Fu-2011}, Fu stated the following fact
\begin{align*}
\sum_{n=0}^{\infty}\mathfrak{B}_k(n)q^n &=\frac{(-q;q)_{\infty}}{(q;q)_{\infty}^{k-1}(-q^k;q^k)_{\infty}}\\
&= \frac{(q^2;q^2)_{\infty}}{(q;q)_{\infty}^k(-q^k;q^k)_{\infty}}\nonumber\\
& \equiv \frac{(q^2;q^2)_{\infty}}{(q^k;q^k)_{\infty}(-q^k;q^k)_{\infty}}\quad (\text{mod}\ p)\nonumber\\
&= \frac{(q^2;q^2)_{\infty}}{(q^{2k};q^{2k})_{\infty}} \quad (\text{mod}\ p).
\end{align*}
Due to Theorem \ref{f-p}, for any prime $p\geq 5$, we have
\begin{equation*}
\sum_{n=0}^{\infty}\mathfrak{B}_k(pn+\frac{p^2-1}{12})q^n \equiv (-1)^{\frac{\pm p-1}{6}}\frac{(q^{2p};q^{2p})_{\infty}}{(q^{2p^{r-1}};q^{2p^{r-1}})_{\infty}} \quad (\text{mod}\ p).
\end{equation*}
That means the lemma holds for $\alpha=1$.
Suppose that lemma holds for $\alpha$. Now we prove the case for $\alpha+1$. For
\begin{equation*}
\sum_{n=0}^{\infty}\mathfrak{B}_k(p^{2\alpha-1}n+\frac{p^{2\alpha}-1}{12})q^n \equiv \left((-1)^{\frac{\pm p-1}{6}}\right)^{\alpha}\frac{(q^{2p};q^{2p})_\infty}{(q^{2p^{r-(2\alpha-1)}};q^{2p^{r-(2\alpha-1)}})_\infty} \quad (\text{mod}\ p).
\end{equation*}
Then
\begin{align}\label{Bp}
\sum_{n=0}^{\infty}\mathfrak{B}_k(p^{2\alpha-1}(pn)+\frac{p^{2\alpha}-1}{12})q^n
&=\sum_{n=0}^{\infty}\mathfrak{B}_k(p^{2\alpha}n+\frac{p^{2\alpha}-1}{12})q^n \nonumber \\
&\equiv \left((-1)^{\frac{\pm
p-1}{6}}\right)^{\alpha}\frac{(q^2;q^2)_\infty}{(q^{2p^{r-2\alpha}};q^{2p^{r-2\alpha}})_\infty}
\quad (\text{mod}\ p).
\end{align}
Using Theorem \ref{f-p} again, we have
\begin{align*}
&\sum_{n=0}^{\infty}\mathfrak{B}_k(p^{2\alpha}(pn+\frac{p^2-1}{12})+\frac{p^{2\alpha}-1}{12})q^n\\
=&\sum_{n=0}^{\infty}\mathfrak{B}_k(p^{2\alpha+1}n+\frac{p^{2\alpha+2}-1}{12})q^n \\
\equiv& \left((-1)^{\frac{\pm
p-1}{6}}\right)^{\alpha+1}\frac{(q^{2p};q^{2p})_\infty}{(q^{2p^{r-(2\alpha+1)}};q^{2p^{r-(2\alpha+1)}})_\infty}
\quad (\text{mod}\ p).
\end{align*}
Therefore the lemma holds for $\alpha+1$.
\qed

According to Lemma \ref{lemma1}, we have the following results.
\begin{thm} For any prime $p \geq 5,n\geq 0$, and $r\geq 1$, if $k=p^r$ is a prime power,
then we have the following two cases:
\begin{enumerate}
\item For $i=1,2,\cdots,p-1$ and $1\leq \alpha\leq r/2$, we have
\begin{equation*}
\mathfrak{B}_k(p^{2\alpha}n+\frac{(12i+p)p^{2\alpha-1}-1}{12}) \equiv 0 \quad (\text{mod}\ p).
\end{equation*}
\item Let $j$ be an integer between $1$ and $p-1$ and $12j+1$ is a quadratic nonresidue
modulo $p$. For $n\geq 0$ and $1\leq \alpha\leq (r-1)/2$, we have
\begin{equation*}
\mathfrak{B}_k(p^{2\alpha+1}n+\frac{(12j+1)p^{2\alpha}-1}{12})\equiv
0\quad (\text{mod}\ p).
\end{equation*}
\end{enumerate}
\end{thm}
\pf According to Lemma \ref{lemma1}, when $1\leq \alpha\leq r/2$, for $i=1,2,\cdots,p-1$, we have
\begin{align*}
\mathfrak{B}_k(p^{2\alpha-1}(pn+i)+\frac{p^{2\alpha}-1}{12})&=\mathfrak{B}_k(p^{2\alpha}n+\frac{(12i+p)p^{2\alpha-1}-1}{12})\equiv
0 \quad (\text{mod}\ p).
\end{align*}

For $1\leq \alpha\leq (r-1)/2$, according to \eqref{Bp} and Theorem
\ref{f-p}, we know that the powers of $q$ modulo $p$ congruent to
 $2\cdot(3k^2+k)/2$ for $-(p-1)/2 \leq k \leq (p-1)/2$ in the
expansion of $(q^2;q^2)_{\infty}$. So we have
\begin{align*}
j &\equiv 2\cdot\frac{3k^2+k}{2}\quad (\text{mod}\ p),\\
12j+1 & \equiv (6k+1)^2\quad (\text{mod}\ p).
\end{align*}
Therefore, if $12j+1$ is a quadratic nonresidue modulo $p$, then we
have
\begin{equation*}
\sum_{n=0}^{\infty}\mathfrak{B}_k(p^{2\alpha}(pn+j)+\frac{p^{2\alpha}-1}{12})q^n
\equiv 0 \quad (\text{mod}\ p).
\end{equation*}\qed

Based on Lemma \ref{lemma1} and the generating
functions of $p(n)$ and $b_{\ell}(n)$, we get the following
congruent relations.

\begin{thm}For any prime $p \geq 5$, $\alpha \geq 1$, and $n \geq 0$, if $k=p^{2\alpha-1}$ is a prime power, then we have
\begin{align}\label{B-b-p}
\sum_{n=0}^{\infty}\mathfrak{B}_k(2p^{2\alpha-1}n+\frac{p^{2\alpha}-1}{12})q^n
&\equiv \left((-1)^{\frac{\pm
p-1}{6}}\right)^{\alpha}\sum_{n=0}^{\infty}b_p(n)q^n \quad
(\text{mod}\ p),\\
\label{Bkpn}
\sum_{n=0}^{\infty}\mathfrak{B}_k(2p^{2\alpha-1}n+\frac{p^{2\alpha}-1}{12})q^n
&\equiv \left((-1)^{\frac{\pm p-1}{6}}\right)^{\alpha}
(q^{p};q^{p})_{\infty}\sum_{n=0}^{\infty}p(n)q^n \quad (\text{mod}\
p).
\end{align}
\end{thm}
\pf Set $r=2\alpha-1$ in Lemma \ref{lemma1}. Then $k=p^{2\alpha-1}$.
So we get \begin{equation*}
\sum_{n=0}^{\infty}\mathfrak{B}_k(p^{2\alpha-1}n+\frac{p^{2\alpha}-1}{12})q^n\equiv
\left((-1)^{\frac{\pm p-1}{6}}\right)^{\alpha}
\frac{(q^{2p};q^{2p})_{\infty}}{(q^2;q^2)_{\infty}} \quad
(\text{mod}\ p).
\end{equation*}
Therefore,
 \begin{align*}
\sum_{n=0}^{\infty}\mathfrak{B}_k(p^{2\alpha-1}(2n)+\frac{p^{2\alpha}-1}{12})q^n
&=\sum_{n=0}^{\infty}\mathfrak{B}_k(2p^{2\alpha-1}n+\frac{p^{2\alpha}-1}{12})q^n\\
&\equiv \left((-1)^{\frac{\pm p-1}{6}}\right)^{\alpha}
\frac{(q^{p};q^{p})_{\infty}}{(q;q)_{\infty}} \quad (\text{mod}\ p).
\end{align*}\qed

Combining \eqref{Bkpn} with the three famous congruences for $p(n)$
given by Ramanujan in \cite{Ramanujan-1919,Ramanujan-2000}
\begin{align}
p(5n+4) &\equiv 0 \quad (\text{mod\ }5),\label{5n+4}\\
p(7n+5) &\equiv 0 \quad (\text{mod\ }7),\label{7n+5}\\
p(11n+6) &\equiv 0 \quad (\text{mod\ }11),\label{11n+6}
\end{align}
we get the following results.

\begin{cor}For $\alpha \geq 1$ and $n \geq 0$, we have
\begin{align*}
\mathfrak{B}_{5^{2\alpha-1}}(2\cdot 5^{2\alpha}n+\frac{101 \cdot 5^{2\alpha-1}-1}{12}) &\equiv 0 \quad (\text{mod}\ 5),\\
\mathfrak{B}_{7^{2\alpha-1}}(2\cdot 7^{2\alpha}n+\frac{127 \cdot 7^{2\alpha-1}-1}{12}) &\equiv 0 \quad (\text{mod}\ 7),\\
\mathfrak{B}_{{11}^{2\alpha-1}}(2\cdot {11}^{2\alpha}n+\frac{155 \cdot {11}^{2\alpha-1}-1}{12}) &\equiv 0 \quad (\text{mod}\ 11).
\end{align*}
\end{cor}
\pf
According to \eqref{Bkpn}, we have
\begin{align*}
\sum_{n=0}^{\infty}\mathfrak{B}_{5^{2\alpha-1}}(2\cdot 5^{2\alpha-1}n+\frac{5^{2\alpha}-1}{12})q^n &\equiv (-1)^{\alpha}(q^5;q^5)_{\infty}\sum_{n=0}^{\infty}p(n)q^n
\quad (\text{mod}\ 5),\\
\sum_{n=0}^{\infty}\mathfrak{B}_{7^{2\alpha-1}}(2\cdot 7^{2\alpha-1}n+\frac{7^{2\alpha}-1}{12})q^n &\equiv (-1)^{\alpha}(q^7;q^7)_{\infty}\sum_{n=0}^{\infty}p(n)q^n
\quad (\text{mod}\ 7),\\
\sum_{n=0}^{\infty}\mathfrak{B}_{{11}^{2\alpha-1}}(2\cdot {11}^{2\alpha-1}n+\frac{{11}^{2\alpha}-1}{12})q^n &\equiv (q^{11};q^{11})_{\infty}\sum_{n=0}^{\infty}p(n)q^n
\quad (\text{mod}\ 11).
\end{align*}
Based on \eqref{5n+4}, \eqref{7n+5}, and \eqref{11n+6}, we get
\begin{align*}
\mathfrak{B}_{5^{2\alpha-1}}(2\cdot 5^{2\alpha-1}(5n+4)+\frac{5^{2\alpha}-1}{12}) &\equiv 0 \quad (\text{mod}\ 5),\\
\mathfrak{B}_{7^{2\alpha-1}}(2\cdot 7^{2\alpha-1}(7n+5)+\frac{7^{2\alpha}-1}{12}) &\equiv 0 \quad (\text{mod}\ 7),\\
\mathfrak{B}_{{11}^{2\alpha-1}}(2\cdot {11}^{2\alpha-1}(11n+6)+\frac{{11}^{2\alpha}-1}{12}) &\equiv 0 \quad (\text{mod}\ 11).
\end{align*}
\qed

Another congruence modulo $p$ for $\mathfrak{B}_k(n)$ can be
directly obtained from Lemma \ref{lemma1}.
\begin{thm}
For any prime $p \geq 5$, $\alpha \geq 1$, and $n \geq 1$, if $k=p^{2\alpha}$ is a prime power, then we have
\begin{equation*}
\mathfrak{B}_k(p^{2\alpha-1}n+\frac{p^{2\alpha}-1}{12}) \equiv 0
\quad (\text{mod}\ p).
\end{equation*}
\end{thm}
\pf Set $r=2\alpha$ in Lemma \ref{lemma1}. Then $k=p^{2\alpha}$. So
we have
\begin{align*}
\sum_{n=0}^{\infty}\mathfrak{B}_k(p^{2\alpha-1}n+\frac{p^{2\alpha}-1}{12})q^n
&\equiv \left((-1)^{\frac{\pm p-1}{6}}\right)^{\alpha} \quad
(\text{mod}\ p).
\end{align*}
\qed

\noindent {\bf Acknowledgements:} This work was supported by the
National Natural Science Foundation of China and the PCSIRT Project
of the Ministry of Education.
%================================================================

\end{document}